\documentclass[12pt]{myiopart}
\usepackage{amsfonts}
\usepackage{amsmath}
\newtheorem{theo}{Theorem}
\newtheorem{prop}{Proposition}{\bf }

\begin{document}
\title{Rigorous results on the threshold network model}

\author{Norio Konno${}^1$, Naoki Masuda${}^{2}$,
Rahul Roy${}^{3}$ and Anish Sarkar${}^{3}$}
\address{
${}^1$ Faculty of Engineering,
Yokohama National University,
79-5, Tokiwadai, Hodogaya, Yokohama, 240-8501, Japan}
\address{${}^2$ Laboratory for Mathematical Neuroscience, RIKEN Brain Science
Institute, 2-1, Hirosawa, Wako, Saitama, 351-0198, Japan}
%
%
\address{
${}^3$ Indian Statistical Institute,
7 SJS Sansanwal Marg, New Delhi 110016, India}


\begin{abstract}
We analyze the threshold network model in which
a pair of vertices with random weights
are connected by an edge when
the summation of the weights exceeds a threshold.
We prove some convergence theorems and central limit theorems
on the vertex degree, degree correlation, and the number of 
prescribed subgraphs. We also generalize some results in 
the spatially extended cases.
\end{abstract}



\maketitle

\section{Introduction}
Recently, real large-scale network data have been analyzed in various
fields, and the corresponding random graphs have been studied. Many 
of these graphs exhibit the power-law
form of the degree distribution, with the power-law exponent typically 
between 2 and 3. In view of this, physicists have proposed
stochastic algorithms for generating growing networks. The underlying 
assumption of these growing networks is that a vertex
with a certain fixed number of edges is added to the graph one by one
at each discrete time step. To obtain a power law, the growth
mechanism is usually supplied by so-called preferential attachment,
which stipulates that each newly introduced edge is more likely to be
connected to a vertex (which already exists in the graph) that has a
larger degree.  For a review of the studies on such models see 
\cite{Albert02,NewmanSIAM}.

However, not all real networks are growing;  algorithms for generating
nongrowing realistic graphs
could be more appropriate for
real situations in which the number of vertices does not change 
rapidly.  In this regard, a type of random graph in which each of 
$n$ vertices is assigned a random variable (weight) was proposed and 
its mean behavior has been analyzed \cite{Caldarelli,Soderberg,Boguna,Masuda_THRESH}.
Interestingly, the power-law degree distribution can emerge even with a weight
distribution that is not equipped with a power-law. For example,
 mean-field analysis suggests a power-law degree distribution
with scaling exponent $-2$ when the vertex
weights are independent
and identically distributed (i.i.d.) random variables obeying
the exponential distribution.

Our first study is regarding such graphs. Formally, our model consists
of $n$ vertices labeled $1,\dotsc, n$ and an i.i.d. sequence of random
variables $X_1, \dotsc, X_n$ with $F$ denoting their common
distribution function. We associate the random variable $X_i$ with the
vertex labeled $i$. Given a fixed threshold value $\theta>0$, we
connect the vertices $i$ and $j$ by an edge $\left<i,j\right>$ if and
only if $X_i + X_j > \theta$ and $i\neq j$.
 Let $G_\theta$ be the random graph thus
produced. A simple coupling argument shows that ${\mathbb P}
\{\left<i,j\right> \in G_{\theta^\prime}\} \leq {\mathbb P}
\{\left<i,j\right> \in G_{\theta}\}$ whenever $\theta^\prime \ge
\theta$.


Let $D_n(i):= \#\{j: \left<i,j\right> \in G_\theta\}$, i.e. the degree
of the vertex $i$ in the graph $G_\theta$. Note $\{D_n(i): i \geq 1\}$
are identical in distribution and let $D_n$ denote a random variable
with this common distribution. The distribution $D_n$ can be obtained
as follows: given $ (n+1) $ vertices, conditioned on the event $ X_1 =
x$, vertices $ j \in \{ 2, \dotsc, n+1 \}$, will connect to the vertex
$ 1 $ if and only if $ X_j > \theta - x $. Therefore, for $ 0 \leq k
\leq n $, we have
\begin{equation}
\label{exdis}
 {\mathbb P} ( D_{n+1} = k ) = \int_{-\infty}^{\infty} {n \choose k} 
[ 1 - F(\theta -x) ]^k [ F (\theta - x)]^{n-k } F(dx).
\end{equation}
Here $F(dx)$ represents the probability measure on the real line
representing the law of $X_1$.
Equation~(\ref{exdis})
allows us to obtain the following asymptotic result:

\begin{theo}
\label{Dn_asymp}
 As $n \to \infty$, 
\begin{equation}
\frac{D_n}{n} \Longrightarrow 1 - F(\theta - X_1).
\end{equation}
\end{theo}
Here $\Longrightarrow$ stands for convergence in distribution.

In data-analysis context, we are often concerned with the degree
correlation between a pair of vertices \cite{Albert02,NewmanSIAM}.
Here we discuss the asymptotic properties of the graph. Our first result is

\begin{theo}
\label{Dn_depend1}
 $\frac{D_n(1)}{n}$ and $\frac{D_n(2)}{n}$ are asymptotically
independent.
\end{theo}

This asymptotic independence breaks under the condition that
the vertices $1$ and $2$ are connected. However, we need to be careful
about this statement, in the sense that we will need some conditions
on the distribution function $F$ of the weight of a vertex. In
particular we assume that

\noindent {\bf Assumption 1: }{\it There exists $u$ and $v$ in the
support of the distribution function $F$
(i.e. for any $\epsilon > 0$, $F(u + \epsilon)
> F(u - \epsilon)$ and $F(v + \epsilon)>F(v - \epsilon)$)
such that $u < \theta/2 < v$ and $u + v > \theta$}\/.
%


\begin{theo}
\label{Dn_depend2}
Given that the vertices $1$ and $2$ are connected, under the above assumption
$\frac{D_n(1)}{n}$ and $\frac{D_n(2)}{n}$ are not asymptotically
independent.
\end{theo}

The importance of the assumption can be understood by looking at the
situation when the assumption does not hold. In that case, the
vertices can be classified into two groups: those which have weights
$\theta/2$ or less and those which have weights more than $\theta/2$. In
the random graph constructed, each of the vertices of the former group
remain isolated, whereas the vertices of the latter group form a
complete graph among themselves. Also, by the strong law of large
numbers, the probability that a vertex belongs to the former group
converges to $F(\theta/2)$ and the probability that a
vertex belongs to the latter group converges to
$1-F(\theta/2)$.  Moreover, given two vertices ($1$ and $2$ (say)) are
connected, they must belong to the latter group and then we have
\begin{eqnarray}
&&{\mathbb P}\left[
\lim_{n\to\infty}\frac{D_n(1)}{n} 
= 1-F(\theta/2),\lim_{n\to\infty} \frac{D_n(2)}{n} = 1-F(\theta/2)
\big| \mbox{vertices }1 \mbox{ and } 2 \mbox{ are connected}\right]\nonumber\\
&&= {\mathbb P}
\left[\lim_{n\to\infty}\frac{D_n(1)}{n} 
= 1-F(\theta/2) \big| \mbox{vertices }1 \mbox{ and } 2
\mbox{ are connected}\right] = 1.
\end{eqnarray}
Thus we obtain conditional asymptotic independence.

A remarkable characteristic of  real graphs is the clustering property 
\cite{Albert02,NewmanSIAM}. The clustering property means an 
abundance of connected triangles in the random graph. The threshold model
exhibits such a clustering property as shown {\it exactly}
in \cite{Boguna,Masuda_THRESH}. 

To formalize this, let $h: {\mathbb R}^3\to{\mathbb R}$ be given by
\begin{equation}
h (x_1,x_2,x_3) :=
1_{ \{ x_1+x_2 > \theta, x_2 + x_3 > \theta, x_3+x_1 > \theta \} }.
\label{hdef-rr}
\end{equation}
Also let 
\begin{eqnarray}
T_n &:=&
 \# \Bigl\{  (i,j,k) : 1 \leq i < j < k \leq n,\Bigr. \nonumber\\
\Bigl. && \;\; X_i + X_j > \theta, 
X_j + X_k > \theta,  X_k + X_i > \theta \Bigr\},\\
F_3 (\theta) &:=& {\mathbb E} \left[ h ( X_1, X_2, X_3) \right]\nonumber\\
& &=  \int_{\mathbb R}  \int_{\mathbb R}  \int_{\mathbb R} 
F(dx_1) F(dx_2) F(dx_3)h(x_1,x_2,x_3)\nonumber\\
&& =  {\mathbb P} ( X_1+X_2 > \theta , X_2+X_3 > \theta , X_3+X_1 > \theta),\\
%
%
\zeta_1 (F) &:=&
{\mathbb E}\left[ \left(\int_{\mathbb R}\int_{\mathbb R}
F(dx_2)F(dx_3) h(X_1, x_2, x_3)\right)^2
\right] - \bigl(F_3 (\theta) \bigr)^2 > 0.
\end{eqnarray}
Note that $T_n$ counts all triangles in the graph with $n$ vertices.

The asymptotic results are

\begin{theo}
\label{triangle_as}
As $ n \to \infty $,
\begin{enumerate}
\item[(a)] $\displaystyle{ \frac{ T_n }{{ n \choose 3}} \to
F_3(\theta) \text{ almost surely; } }$

\item[(b)] $ \displaystyle{\sqrt{n} \Bigl[ \frac{T_n}{{ n \choose 3}} - 
F_3 (\theta) \Bigr]  \Rightarrow \sqrt{3 \zeta_1 (F)} Z}$
where $Z$ is a standard normal random variable.
\end{enumerate}
\end{theo}
The method used to show the above results may be generalized easily to 
obtain a count of not only triangles but
subgraphs in $G_\theta$ isomorphic to a fixed graph. 
In data-analysis contexts, a fixed subgraph is 
called a motif of the graph.
Depending on the types of real networks
(e.g. Internet, gene networks, neural networks, social networks),
there are some sorts of small motifs
that appear in an entire graph significantly more than in the random graphs.
These motifs are relevant to
functional roles such as signal transduction and information processing
apposite to each application
\cite{Itzkovitz,Milo02}.
Our results, Theorems \ref{VE_as} and \ref{VE_clt} in Section 4,
obtain limit theorems
for the motifs of the threshold model.

Besides the extension to general graphs, the above theorem may be
extended for local triangle counts, i.e., the number of triangles
containing a specified vertex. The clustering coefficient, which is a
quantity often used to evaluate the degree of clustering property, is
defined as the normalized number of locally counted triangles averaged
over all the vertices. Our local results below show that the vertex-wise
clustering coefficient satisfies the central limit theorem, which
validates the use of expectation in physics community.

We fix the vertex  $ 1 $ and define $ T_n (1) $ as the number of triangles in 
the set of vertices $ \{  1,2, \dotsc, n, n+1 \} $ of which the vertex $ 1 $ is a site. 
In other words, 
\begin{equation}
\label{def:loc_tringles}
T_n (1) := \# \{ 2 \leq i < j \leq n+1 : h(X_1, X_i, X_j)=1\}.
\end{equation}

\begin{theo}
\label{thm:local_triangles}
As $ n \to \infty $, 
\begin{equation}
\frac{ T_n (1)}{ {n\choose2}} \Rightarrow \int_{\mathbb R} 
\int_{\mathbb R} F(dx_1) F(dx_2) h(X,x_1,x_2)
\end{equation}
where $X$ is an independent random variable identical in distribution
to $X_1$.
 \end{theo}

Besides the above dimensionless random graph model, a spatial model has 
been proposed in \cite{Masuda_SPTH}. Consider a homogeneous Poisson point 
process of intensity $ \lambda $ on ${\mathbb R}^d$. 
We insist that the origin, ${\mathbf 0}$, is a point of the process. Let 
$ \{ {\mathbf 0} = \xi_0, \xi_1, \xi_2, \dotsc \} $ be an enumeration of the 
points of the process.  Associated with each point $ \xi_i $ is a random variable
$ X_i $. We assume that $ \{ X_0, X_1, X_2, \dotsc \} $ is an
i.i.d. sequence of random variables with common 
distribution function $ F $. The random graph $G_{\theta, \beta}$ is obtained 
by connecting  $ \xi_i $ and $ \xi_j $ by an edge if and only if $ (X_i + X_j) > 
\theta | \xi_i - \xi_j |^{\beta } $ where $ \theta $ and $ \beta $,
the parameters of the model, are real numbers.

Define the degree of the origin in a sphere of radius $r$ as 
\begin{equation}
\label{def:deg0n}
\Delta_r : = \# \{ i \geq 1 : (X_0 + X_i) > \theta | \xi_i  |^{\beta } \text {and } |\xi_i| \leq r \} .
\end{equation}
Given a fixed 
$ x \in {\mathbb R} $, let
\begin{equation}
\label{r-fdef}
f (r;x) = f(|r|;x): = 1 - F ( \theta  |r|^{\beta }  - x ).
\end{equation}
Define 
\begin{equation}
C_r (x) := \int_0^r \tilde{r}^{d-1} f(\tilde{r};x) d\tilde{r},
\end{equation}
and consider the two following cases:
\begin{enumerate}
\item As $ r \to \infty $, $ C_r (x) \to C(x) := \int_{0}^{\infty}
\tilde{r}^{d-1} f(\tilde{r},x)d\tilde{r} < \infty $ for every $ x \in {\mathbb R} $.
\item There exists a sequence $ \{C_r\}$ and a function $g(x)$ such
that $ C_r \to \infty $ and $ \frac{ C_r (x) - C_r }{ \sqrt{C_r}}
\to g(x) $ as $ r \to \infty $ for every $ x \in {\mathbb R} $.
\end{enumerate}

For the first case, we have
\begin{theo}
\label{finitecase}
If $ C_r (x) \to C(x) < \infty $ for every $ x \in {\mathbb R} $ as $
r \to \infty $, then we have
\begin{equation}
\Delta_r \Rightarrow \Delta
\end{equation}
where the characteristic function of the random variable $\Delta$ is
given by $ \phi_{\Delta} (t) = \int_{\mathbb R} F(dx) \exp \bigl( -
\lambda c_d C(x) ( 1 - \exp (it) ) \bigr) $ where $ c_d $ represents
the volume of the $(d-1)$-dimensional unit sphere.
\end{theo}
\noindent{\bf Remark } The random variable $\Delta$ represents the degree 
of the origin in the random graph. 

For the second case, we have
\begin{theo}
\label{clt}
Suppose that there exists a sequence  $ \{C_r\}$
such that $ C_r \to \infty $ and 
$ \frac{ C_r (x) - C_r}{ \sqrt{ C_r }}  \to g(x) $ as $r \to \infty $ for 
every $ x \in  {\mathbb R} $. We have
\begin{equation}
\frac{ \Delta_r - \lambda c_d C_r }{ \sqrt{ \lambda c_d C_r} } \Rightarrow  
Z + \sqrt{\lambda c_d} g(X_0) \text{ as } r \to \infty
\end{equation}
where
$ c_d $ is as defined in Theorem \ref{finitecase}. 
\end{theo}
We end this section with an example of $ F $ satisfying the condition 
in Theorem \ref{clt}. Fix $ \beta = 1 $ and $ d = 2 $. 
Define $ F : [0, \infty) \to [0,1] $ by
\begin{equation}
F(x) := 1 - C x^{ - \alpha }
\end{equation}
where $ 0 < \alpha < 2 $ and $ C > 0 $, and
\begin{equation}
C_r := C \theta^{-\alpha} r^{ 2 - \alpha} / (2 - \alpha).
\end{equation}
Simple computations can be carried out to verify that the conditions
of Theorem \ref{clt} are satisfied in this case with $g(x) = 0$ for all 
$x\in \mathbb R$. In the next few sections we prove our results.

\section{Degree $D_n$ of a vertex}
To prove Theorem \ref{Dn_asymp} observe that, since $ (n+1)/n \to 1 $, 
it is enough to show that $ D_{n+1}/n $ converges weakly to the required 
random variable. For $\xi \in \mathbb R$, we have from
equation~(\ref{exdis})
\begin{eqnarray}
\lefteqn{ {\mathbb E}\left[\exp \left(\tfrac{it D_{n+1}}{n}\right)\right] }\nonumber\\
& = &  \sum_{k=0}^n  \exp \left(\frac{itk}{n}\right)
\int_{-\infty}^\infty 
{n \choose k} [1-F(\theta - x)]^k  [F(\theta - x)]^{n-k} F(dx)\nonumber\\
& = &  \int_{-\infty}^\infty \sum_{k=0}^n {n \choose k} 
\exp\left(\frac{itk}{n}\right) [1-F(\theta - x)]^k 
[F(\theta - x)]^{n-k} F(dx)\nonumber\\
& = & \int_{-\infty}^\infty \left[(1- F(\theta - x)) 
\exp\left(\frac{it}{n}\right)
 + F(\theta - x) \right]^n F(dx)\nonumber\\
& = & \int_{-\infty}^\infty \left[(1- F(\theta - x))
\left( 1 +  \frac{i t}{n} 
+ o\left(\frac{1}{n}\right)\right) + F(\theta - x) \right]^n F(dx)\nonumber\\
& = & \int_{-\infty}^\infty \left[ 1 +  \frac{i t}{n}(1- F(\theta - x)) 
+ o\left(\frac{1}{n}\right)\right]^n F(dx)\nonumber\\
& \rightarrow & \int_{-\infty}^\infty \exp({i t (1- F(\theta - x))})
F(dx) 
\text{ as } n \to \infty,
\end{eqnarray}
where the limit follows from dominated convergence theorem because
\begin{eqnarray}
\lefteqn{\left|(1- F(\theta - x)) \exp\left(\frac{it}{n}\right) + F(\theta - x)\right|^2}\nonumber\\
& = &  [(1- F(\theta - x)) \cos(t/n) + F(\theta - x)]^2 + (1- F(\theta - x))^2 
\sin^2(t/n)\nonumber\\
& = & (1- F(\theta - x))^2 + 2 F(\theta - x)(1- F(\theta - x))\cos(t/n) 
+ (F(\theta - x))^2\nonumber\\
& \leq & 1.
\end{eqnarray}
This proves Theorem \ref{Dn_asymp}. 

\section{Pair correlation}
Now suppose there are $n+2$ vertices labeled $1,2,\dotsc n, n+1, n+2$ 
with corresponding  random variables $X_1, \dotsc , X_{n+2}$. Define
\begin{eqnarray}
D_{n,1} &:=& D_{n+2}(n+1) = \#\{j: 1 \leq j \leq n+2, j \neq n+1, X_j 
+ X_{n+1} >  \theta\}\\
D_{n,2} &:=& D_{n+2}(n+2) =  \#\{j: 1 \leq j \leq n+2, j \neq n+2, 
X_j + X_{n+2} > \theta\}.
\end{eqnarray}

We have, for $0 \leq k, l \leq n+1$,
\begin{equation}
{\mathbb P}\left[D_{n,1} = k, D_{n,2} = l\right] = {\mathbb P} 
\left[D_{n,1} = l, D_{n,2} = k\right].
\end{equation}

For $ k > l$
\begin{eqnarray}
\lefteqn{{\mathbb P} \left[D_{n,1} = k, D_{n,2} = l\right]}\nonumber\\
& = & \int_{-\infty}^\infty \int_{-\infty}^\infty 1_{\{a > b\}} F(da) F(db) 
{\mathbb P} \left[\#\{1 \leq j \leq n : a+ X_j > \theta\} + 
1_{\{a + b > \theta\}} =k,\right.\nonumber\\
&  & \left.
 \qquad \qquad \qquad \qquad \qquad \qquad \quad \#\{1 \leq j \leq n 
: b+ X_j > \theta\} + 1_{\{a + b > \theta\}} =l \right]\nonumber\\
& = &  \int_{-\infty}^\infty \int_{-\infty}^\infty 1_{\{a > b\}}  
1_{\{a+b > \theta\}} F(da) F(db) {\mathbb P} \left[ \#\{1 \leq j \leq n 
: a+ X_j > \theta\} =k - 1,\right.\nonumber\\
& & \left. \qquad \qquad \qquad \qquad \qquad \qquad \quad \# \{1 \leq j \leq n 
: b+ X_j > \theta\} =l -1 \right]\nonumber\\
& & + \int_{-\infty}^\infty \int_{-\infty}^\infty 1_{\{a > b\}}  
1_{\{a+b \le \theta\}} F(da) F(db) {\mathbb P} \left[ \#\{1 \leq j \leq n 
: a+ X_j > \theta\} =k,\right.\nonumber\\
& & \left.\qquad \qquad \qquad \qquad \qquad \qquad \qquad \quad \#\{1 
\leq j \leq n : b+ X_j > \theta\} =l \right]\nonumber\\ 
& = & \int_{-\infty}^\infty \int_{-\infty}^\infty 1_{\{a > b\}}  
1_{\{a+b > \theta\}} g_n(\theta; a,b; k-1, l-1) F(da) F(db)\nonumber\\
& & +  \int_{-\infty}^\infty \int_{-\infty}^\infty 1_{\{a > b\}}  
1_{\{a+b \le \theta\}} g_n(\theta; a,b; k, l)F(da) F(db),
\end{eqnarray}
where, for $a \ge b$,  
\begin{align}
\lefteqn{ g_n(\theta; a,b; k, l)} \nonumber\\
 := & \medspace {\mathbb P} \left[
\#\{1 \leq j \leq n : a+ X_j > \theta\} =k, 
\#\{1 \leq j \leq n : b+ X_j > \theta\} =l \right] \nonumber\\
 = &  \medspace {\mathbb P} \left[
\#\{1 \leq j \leq n :  X_j > \theta - a\} =k, 
\#\{1 \leq j \leq n :  X_j > \theta - b\} =l \right] \nonumber\\
 = &  \medspace {\mathbb P} \left[ \#\{1 \leq j \leq n :  X_j > \theta - b\} =l, 
\#\{1 \leq j \leq n : \theta - a < X_j \leq \theta - b\} =k -l \right] \nonumber\\
 = &  \medspace 
\left\{\begin{array}{ll}
\frac{n!}{l! (k-l)! (n-k)!}
\left[1- F(\theta -b)\right]^l  \left[F(\theta - b) 
- F(\theta -a)\right]^{k-l} \left[F(\theta - a) \right]^{n-k}
& \mbox{ if } k\ge l\ge 0\\
0 & \mbox{ otherwise}. 
\end{array}\right.
\end{align}

For $k = l$ similar calculations show
\begin{eqnarray}
\lefteqn{{\mathbb P} \left[D_{n,1} = k, D_{n,2} = k\right]}\nonumber\\
& = &
\int_{-\infty}^\infty \int_{-\infty}^\infty 1_{\{a > b\}}  1_{\{a+b > \theta\}}
g_n(\theta; a,b; k-1, k-1) F(da) F(db)\nonumber\\
& & +  \int_{-\infty}^\infty \int_{-\infty}^\infty 1_{\{a > b\}}  1_{\{a+b \le \theta\}} 
g_n(\theta; a,b; k, k)F(da) F(db)\nonumber\\
&  & +\int_{-\infty}^\infty \int_{-\infty}^\infty 1_{\{a \le b\}}  1_{\{a+b > \theta\}}
g_n(\theta; b,a; k-1, k-1) F(da) F(db)\nonumber\\
& & +  \int_{-\infty}^\infty \int_{-\infty}^\infty 1_{\{a \le b\}}  
1_{\{a+b \le \theta\}} g_n(\theta; b,a; k, k)F(da) F(db).
\end{eqnarray}

Now let
\begin{equation}
d_n(1) := \#\{j : 1 \leq j \leq n, X_j + X_{n+1} > \theta\}
\end{equation}
and
\begin{equation}
d_n(2) := \#\{j : 1 \leq j \leq n, X_j + X_{n+2} > \theta\};
\end{equation}
then
\begin{equation}
D_{n,1} \geq d_n(1) \geq D_{n,1} - 1
\end{equation}
and
\begin{equation}
D_{n,2} \geq d_n(2) \geq D_{n,2} - 1.
\end{equation}
Thus the asymptotic distributions of $\left(\tfrac{D_{n,1}}{n},
\tfrac{D_{n,2}}{n}\right)$ and $\left(\tfrac{d_{n}(1)}{n},
\tfrac{d_{n}(2)}{n}\right)$ are identical, and we work out pair
correlation with $d_n(1)$ and $d_n(2)$ instead of $D_{n,1}$ and
$D_{n,2}$.

Now, for all $0 \leq k, l \leq n$,
\begin{equation}
{\mathbb P} \left[d_n(1)=k, d_n(2)=l\right]
= {\mathbb P} \left[d_n(1)=l, d_n(2)=k\right],
\end{equation}
and, for $k > l$
\begin{equation}
{\mathbb P}\left[d_n(1)=k, d_n(2)=l\right] = 
\int_{-\infty}^\infty \int_{-\infty}^\infty 1_{\{a > b\}}   g_n(\theta; a,b; k, l) F(da) F(db),
\end{equation}
while, for $k = l$
\begin{eqnarray}
{\mathbb P}\left[d_n(1)=k, d_n(2)=k\right]
& = &  \int_{-\infty}^\infty \int_{-\infty}^\infty 1_{\{a > b\}}   
g_n(\theta; a,b; k, k)F(da) F(db)\nonumber\\
& & + \int_{-\infty}^\infty \int_{-\infty}^\infty 1_{\{a \le b\} }   
g_n(\theta; b,a; k, k)F(da) F(db)\nonumber\\
& = & 2\int_{-\infty}^\infty \int_{-\infty}^\infty 1_{ \{a > b\}}   
g_n(\theta; a,b; k, k)F(da) F(db)\nonumber\\
&& + \int_{-\infty}^\infty \int_{-\infty}^\infty 1_{ \{a = b\}}   
g_n(\theta; a,b; k, k)F(da) F(db).
\end{eqnarray}

Now let $ -\infty < s, t < \infty$. To derive pair independence, 
let us consider the characteristic function. We have
\begin{equation}
{\mathbb E}\left[\exp ({isd_n(1) + itd_n(2)}) \right] =\sum_{k=0}^n 
\sum_{l=0}^n \exp( {isk + itl}){\mathbb P}\left[d_n(1)=k,
  d_n(2)=l\right].
\end{equation}
We break the above double sum into three parts, when (i) $k > l$, (ii) 
$k < l$ and (iii) $k = l$.

For (i) we have
\begin{eqnarray}
\lefteqn{\sum_{k,l : 0 \leq l < k \leq n}  \exp( {isk + itl}) {\mathbb P} 
\left[d_n(1)=k, d_n(2)=l\right]}\nonumber\\
& = & \sum_{k,l : 0 \leq l < k \leq n}  \exp ({isk + itl} )\int_{-\infty}^\infty 
\int_{-\infty}^\infty 1_{\{a > b\}}   g_n(\theta; a,b; k, l)F(da) F(db)\nonumber\\
& = &  \int_{-\infty}^\infty \int_{-\infty}^\infty 1_{\{a > b\}} F(da) F(db)
\sum_{k,l : 0 \leq l < k \leq n}  \exp( {isk + itl}) g_n(\theta; a,b; k, l).
\end{eqnarray}
The term 
\begin{eqnarray}
\lefteqn{\sum_{k,l : 0 \leq l < k \leq n} \exp ( {isk + itl}) g_n(\theta; a,b; k, l)}\nonumber\\
&=& \sum_{k,l : 0 \leq l < k \leq n}  \exp( {isk + itl}) \frac{n!}{l! (k-l)! (n-k)!} \nonumber\\
& & \qquad \qquad \times \Bigl[1- F(\theta -b) \Bigr]^l \Bigl[F(\theta - b) - F(\theta -a) 
\Bigr]^{k-l}  \Bigl[F(\theta - a) \Bigr]^{n-k}\nonumber\\
%
%
&=& \sum_{k,l : 0 \leq l \leq k \leq n}  \exp ({is(k-l) + il(s+t)} )
\frac{n!}{l! (k-l)! (n-k)!} \nonumber\\
& & \qquad \qquad \times \Bigl[1- F(\theta -b) \Bigr]^l \Bigl[F(\theta - b) - F(\theta -a) 
\Bigr]^{k-l}  \Bigl[F(\theta - a) \Bigr]^{n-k}\nonumber\\
& & - \sum_{k,l : 0 \leq l = k \leq n}  \exp ({is(k-l) + il(s+t)})
\frac{n!}{l! (k-l)! (n-k)!} \nonumber\\
& & \qquad \qquad \times \Bigl[1- F(\theta -b) \Bigr]^l \Bigl[F(\theta - b) - F(\theta -a) 
\Bigr]^{k-l}  \Bigl[F(\theta - a) \Bigr]^{n-k}\nonumber\\
& = & \Bigl[ \exp ({i(s+t)}) \bigl(1-F(\theta - b) \bigr) + \exp ({is}) \bigl(F(\theta - b) 
- F(\theta -a) \bigr) + F(\theta - a) \Bigr]^n\nonumber\\
& & \qquad \qquad - \sum_{k=0}^n  \exp ({ ik(s+t)})
\frac{n!}{k!  (n-k)!} [1- F(\theta -b)]^k [F(\theta - a)]^{n-k}\nonumber\\
& = & \Bigl[ \exp ({i(s+t)}) \bigl(1-F(\theta - b) \bigr) + \exp(is) \bigl(F(\theta - b) 
- F(\theta -a) \bigr) + F(\theta - a) \Bigr]^n\nonumber\\
& & \qquad \qquad - \Bigl[ \exp ({i(s+t)}) \bigl(1-F(\theta - b) \bigr)+ 
F(\theta - a) \Bigr]^n.
\end{eqnarray}
Thus
\begin{eqnarray}
\lefteqn{\sum_{k,l : 0 \leq l < k \leq n}  \exp ({isk + itl}) 
{\mathbb P}\left[d_n(1)=k, d_n(2)=l\right]}\nonumber\\
&=& \int_{-\infty}^\infty \int_{-\infty}^\infty 1_{\{a > b\}}F(da) F(db)\nonumber\\
& & \Bigl(\bigl[\exp({i(s+t)}) (1-F(\theta - b)) + \exp(is)
(F(\theta - b) - F(\theta -a)) 
+ F(\theta - a)\bigr]^n\nonumber\\
& & \qquad \quad - \bigl[\exp( {i(s+t)}) (1-F(\theta - b))+ F(\theta - a) \bigr]^n\Bigr).
\end{eqnarray}
Similarly, for (ii),
\begin{eqnarray}
\lefteqn{\sum_{k,l : 0 \leq k < l \leq n}  \exp ({isk + itl}) {\mathbb P} 
\left[d_n(1)=k, d_n(2)=l\right]}\nonumber\\
&=& \int_{-\infty}^\infty \int_{-\infty}^\infty 1_{\{ a > b \}}F(da) F(db)\nonumber\\
& & \Bigl( \bigl[ \exp({i(s+t)}) (1-F(\theta - b))  + \exp ({it}) (F(\theta - b) - 
F(\theta -a)) + F(\theta - a) \bigr]^n\nonumber\\
& & \qquad \quad - \bigl[ \exp ({i(s+t)}) (1-F(\theta - b))+ F(\theta - a) \bigr]^n\Bigr).
\end{eqnarray}
For (iii),
\begin{eqnarray}
\lefteqn{\sum_{k,l : 0 \leq k = l \leq n}  \exp ({isk + itk}) {\mathbb P} 
\left[d_n(1)=k, d_n(2)=k\right]}\nonumber\\
&=& 2\int_{-\infty}^\infty \int_{-\infty}^\infty 1_{\{ a > b \}}F(da) F(db)\nonumber\\
& & \qquad \Bigl(\sum_{k=0}^n  \exp ({isk + itk}) \frac{n!}{k! (n-k)!}(1-F(\theta - b))^k 
(F(\theta - a))^{n-k}\Bigr) \nonumber\\
&&+ \int_{-\infty}^\infty \int_{-\infty}^\infty 1_{\{ a = b \}}F(da) F(db)\nonumber\\
& & \qquad \Bigl(\sum_{k=0}^n  \exp ({isk + itk}) \frac{n!}{k! (n-k)!}(1-F(\theta - b))^k 
(F(\theta - a))^{n-k}\Bigr)\nonumber\\
& = & 2\int_{-\infty}^\infty \int_{-\infty}^\infty 1_{\{ a > b \}}F(da) F(db)
\bigl[(1-F(\theta - b)) \exp( {i(s+t)}) + F(\theta - a)\bigr]^{n}\nonumber\\
&&+ \int_{-\infty}^\infty \int_{-\infty}^\infty 1_{\{ a = b \}}F(da) F(db)
\bigl[(1-F(\theta - b)) \exp( {i(s+t)}) + F(\theta - a)\bigr]^{n}.
\end{eqnarray}
Combining the above and
%
%
taking $s = u/n$ and $t = v/n$, we have
\begin{eqnarray}
\lefteqn{{\mathbb E}\left[ \exp ({iu(d_n(1)/n) + iv(d_n(2)/n)}) \right]}\nonumber\\
&=& 
\int_{-\infty}^\infty \int_{-\infty}^\infty 1_{\{ a > b \}}F(da)
F(db)\nonumber\\
&&\left\{
\left[\exp\left(\tfrac{i(u+v)}{n}\right) \left(1-F(\theta - b) \right) + \exp\left(\tfrac{iu}{n}\right) 
\bigl( F(\theta - b)- F(\theta - a) \bigr) + F(\theta - a) \right]^n
\right.
\\
& & +\left. \left[\exp\left(\tfrac{i(u+v)}{n}\right) \left(1-F(\theta - b) \right) + \exp\left(\tfrac{iv}{n}\right) 
\bigl( F(\theta - b)- F(\theta - a) \bigr) + F(\theta - a) \right]^n
\right\}\nonumber\\
&& +\int_{-\infty}^\infty \int_{-\infty}^\infty 1_{\{ a = b \}}F(da)
F(db)
\left[\exp\left(\tfrac{i(u+v)}{n}\right) \left(1-F(\theta - b) \right)
+ F(\theta - a) \right]^n\nonumber\\
&=& \int_{-\infty}^\infty \int_{-\infty}^\infty 1_{\{ a > b \}}F(da) F(db)
\left\{
\left[ \left(1 + \frac{i(u+v)}{n} + \frac{O(1)}{n^2} \right) 
\left(1-F(\theta - b) \right) \right.\right. \nonumber\\
& & \qquad \qquad + \left.\left. \left(1 + \frac{iu}{n} + \frac{O(1)}{n^2} \right) 
\left(F(\theta - b)- F(\theta - a) \right) + F(\theta - a)\right]^n \right.\nonumber\\
&&  + \left[ \left(1 + \frac{i(u+v)}{n} + \frac{O(1)}{n^2} \right) 
\left(1-F(\theta - b) \right) \right.\nonumber\\
&& \left.\left. +\left(1 + \frac{iv}{n} + \frac{O(1)}{n^2} \right) 
\left(F(\theta - b)- F(\theta - a) \right)
+ F(\theta -
a)\right]^n\right\}\nonumber\\
&&
+\int_{-\infty}^\infty \int_{-\infty}^\infty 1_{\{a=b\}}F(da) F(db)\left[ \left(1 + \frac{i(u+v)}{n} + \frac{O(1)}{n^2} \right) 
\left(1-F(\theta - b) \right) +
F(\theta - a)\right]^n.
%
%
\end{eqnarray}
Letting $n \to \infty$, we see that
\begin{eqnarray}
\label{lim:joint}
&&\lefteqn{\lim_{n \to \infty}{\mathbb E}\left[
\exp( {iu(d_n(1)/n) + iv(d_n(2)/n)}) 
\right]} \nonumber \\
& = & 
\int_{-\infty}^\infty \int_{-\infty}^\infty 1_{\{ a > b \}}F(da) F(db)
\left\{
\exp \left[{i(u+v) (1-F(\theta -b)) + iu (F(\theta-b)-F(\theta-a))}\right]
\right.\nonumber \\
&&\left.
+\exp \left[{i(u+v) (1-F(\theta-b)) + iv (F(\theta-b)-F(\theta-a))}\right]
\right\}\nonumber \\
&&+ \int_{-\infty}^\infty \int_{-\infty}^\infty 1_{\{a=b\}}F(da)F(db)
\exp \left[{i(u+v) (1-F(\theta-b))}\right]
\nonumber\\
&=&  \int_{-\infty}^\infty \exp ({iv (1 - F(\theta -b))}) F(db)
\int_{-\infty}^\infty \exp ({iu (1 - F(\theta -a))}) F(da).
\end{eqnarray}
Thus combining with the result in Section 2 we see that
\begin{eqnarray}
\label{equality:joint}
\lefteqn{ \lim_{n \to \infty}{\mathbb E} \left[\exp ( {iu(d_n(1)/n) + iv(d_n(2)/n)}) 
\right] } \nonumber \\
& = &  {\mathbb E} \left[\exp (iu (1 - F(\theta - X_1))) \right]
{\mathbb E} \left[\exp (iv (1 - F(\theta - X_1))) \right] \nonumber \\
& = &
\lim_{n \to \infty} {\mathbb E}\left[ \exp({iu(D_n(1)/n)})\right]
\lim_{n \to \infty} {\mathbb E}\left[ \exp({iv(D_n(2)/n)})\right], 
\end{eqnarray}
i.e. we obtain the asymptotic independence as enunciated in Theorem
\ref{Dn_depend1}.

Now we obtain Theorem \ref{Dn_depend2}.
The joint conditional probability  distribution of the weights of two fixed 
vertices, provided that these vertices are connected by an edge,
is given by
\begin{eqnarray}
\label{cond_distrn:joint}
H(da, db)
& = & \frac{ 1  (a+b> \theta) }{ \alpha_F (\theta) } F(da) F(db)
\end{eqnarray}
where $ \alpha_F (\theta) := \int_{-\infty}^{\infty}
\int_{-\infty}^{\infty} 1 (\tilde{a}+\tilde{b}> \theta) F(d\tilde{a})
F(d\tilde{b})$ is the normalizing constant. It is actually the
probability that two vertices share an edge. Consequently, the
conditional probability distribution of weight of a vertex
provided that it shares an edge with another vertex is given by
\begin{eqnarray}
\label{cond_distrn:marginal}
G(da) & = & \int_{-\infty}^{\infty} H(da, db) \nonumber \\
& = & F(da) \frac{ \int_{-\infty}^{\infty} 1(a+b>\theta) F(db)}
{ \alpha_F (\theta) } \nonumber  \\
& = & \frac{1 - F(\theta -a)}{ \alpha_F (\theta) } F(da) .
\end{eqnarray}

In order to study the asymptotic dependence of $ d_n(1)/n $ and $ d_n(2)/n $,
we need to consider the difference between the expressions
in the left and right sides of equation~(\ref{lim:joint}). Using equations
(\ref{cond_distrn:joint}) and (\ref{cond_distrn:marginal}), 
the asymptotic limit of this difference is
\begin{eqnarray}
\label{cond:asym_ind}
&&\lim_{n\to\infty}{\mathbb E}\left[\exp
 ({iu(d_n(1)/n)+iv(d_n(2)/n)})\right]
 \nonumber \\
& & \quad - \lim_{n\to\infty}{\mathbb E}\left[\exp ({iu(D_n(1)/n)})\right]
\lim_{n\to\infty}{\mathbb E}\left[ \exp ({iv(D_n(2)/n)})\right]\nonumber\\
&=& \int^{\infty}_{-\infty} \left( H(da, db) - G(da)G(db)\right)
e^{iu(1-F(\theta-a))} e^{iv(1-F(\theta-b))}.
\end{eqnarray}
If the resulting conditional random variables are asymptotically independent, then 
for all $ u,v \in {\mathbb R} $, the right hand 
side of equation~(\ref{cond:asym_ind}) must vanish. Now we claim that 
this will imply that the probability
 measures $ H (da, db) $ and $ G(da) G(db) $ 
on $ {\mathbb R}^2 $ are identical. Indeed, let $ (M_1,M_2) $ and 
$ (R_1, R_2) $ be two random vectors on $ {\mathbb R}^2 $
whose distributions are given by the probability measures $ H (da, db) $ 
and $ G(da) G(db) $, respectively. Consider the map $\tilde{f}$ from 
$ {\mathbb R}^2 \to {\mathbb R}^2$ defined by 
$\tilde{f}: (a,b) \to (1 - F(\theta -a), 1 -F(\theta - b) ) $. Then,
the supposition that equation~(\ref{cond:asym_ind}) vanishes
implies that the characteristic function 
of $\tilde{f}(M_1,M_2) $ is same as that of $\tilde{f}(R_1, R_2)$. Hence their distributions
are also same. Therefore, 
\begin{eqnarray}
{\mathbb P} \lefteqn{ \left[M_1 \leq \beta_1, M_2 \leq
\beta_2\right]}\nonumber\\
& = & {\mathbb P}\left[1 - F( \theta - M_1)
\leq 1 - F(\theta - \beta_1) , 1 - F( \theta - M_2) \leq 1 - F(\theta
- \beta_2)\right]\nonumber\\
& = & {\mathbb P} \left[1 - F( \theta -
R_1) \leq 1 - F(\theta - \beta_1) , 1 - F( \theta - R_2) \leq 1 -
F(\theta - \beta_2)\right]\nonumber\\
& = & {\mathbb P}\left[R_1 \leq
\beta_1, R_2 \leq \beta_2 \right];
\end{eqnarray}
{\it i.e. } $H(da, db)$ and $G(da)G(db)$ are identical.

Now, we claim that, because of
Assumption 1 (of Section 1), the probability
measures $ H (da, db) $ and $ G(da) G(db) $ cannot be the same.
Indeed, if the measures $ H (da, db) $ and $ G(da) G(db) $ are same,
then for any subset $ A \subseteq {\mathbb R}^2 $, $ \int_A H(da, db)
= \int_A G(da) G(db) $.  Using the equations (\ref{cond_distrn:joint})
and (\ref{cond_distrn:marginal}), we have
\begin{equation}
\int_A  \left[ \frac{ 1  (a+b> \theta) }{ \alpha_F (\theta) } - 
\frac{ 1 - F(\theta -a)  }{ \alpha_F (\theta) } \frac{ 1 - F(\theta -b) }
{ \alpha_F (\theta) }  \right] F(da) F(db) = 0.
\end{equation}
This will imply that 
\begin{equation}
1  (a+b> \theta)  \alpha_F (\theta)
 -  (1 - F(\theta -a) )  (1 - F(\theta -b) ) 
= 0
\end{equation}
$ F \times F $ almost surely.

Now fix $u$ and $v$ as in Assumption 1 and let $\epsilon$ be such that
$u+ v - 4\epsilon > \theta$ and $u+4 \epsilon < \theta/2 < v - 4
\epsilon$. Let $a, b \in (u-\epsilon, u+\epsilon]$. Then
\begin{equation}
1(a+b > \theta)  \alpha_F (\theta)= 0.
\end{equation}
However, $\theta - a <$ $\theta - u + \epsilon <$ $v - 3 \epsilon$
$<v-\epsilon$,
so $ 1 - F(\theta -a)$  $\geq 1 - F(v-\epsilon)$
$\geq F(v+\epsilon) - F(v-\epsilon) > 0$.
Similarly, $ 1 - F(\theta -b) > 0$. Thus
\begin{equation}
(1 - F(\theta -a) )  (1 - F(\theta -b) ) > 0.
\end{equation}
Hence, on a set of probability at least
 $(F(v+\epsilon) - F(v-\epsilon))^2>0$ we have
\begin{equation}
1  (a+b> \theta)  \alpha_F (\theta) 
- (1 - F(\theta -a) )  (1 - F(\theta -b) ) < 0.
\end{equation}
This contradiction completes the proof.

\section{Triangles} 

Now we study the number of triangles in the graph. The number of
triangles $T_n$ can be represented as a U-statistic (see
\cite{Serfling}).  The kernel function $ h : {\mathbb R}^3 \to
{\mathbb R} $ defined in equation (\ref{hdef-rr}) is clearly a
symmetric function of $x_1$, $x_2$, $x_3$. Then, we have
\begin{equation}
T_n = {n \choose 3} \times \frac{1}{{n \choose 3}} \sum_{ (i,j,k)
\in C } h (X_i,X_j,X_k) = {n \choose 3} U_n
\end{equation}
where $C :=
\{(i,j,k) : 1\le i<j<k\le n\}$ is the
collection of all possible triplets and $U_n :=
\frac{1}{{n \choose 3}} \sum_{ (i,j,k)\in C } h (X_i,X_j,X_k)$
is the U-statistic
obtained from the kernel $ h $.  Theorem \ref{triangle_as}(a) can be
easily derived from Theorem A of \cite{Serfling} (1980), page 190.

Next define $ h_1 : {\mathbb R} \to {\mathbb R}$ as follows:
\begin{eqnarray}
\label{assocfn}
 h_1 (x) & := & {\mathbb E}\left[ h(x, X_2, X_3) \right] \nonumber\\
& = & \mathbb P \left[
 X_2 > \theta - x, X_3 > \theta - x, X_2 +X_3 > \theta \right] \nonumber\\
& = &  \mathbb P \left[
\min \{ X_2, X_3 \}  > \theta - x, X_2 +X_3 > \theta \right]. 
\end{eqnarray}
Noting that
\begin{equation}
{\mathbb E}\left[h_1 (X_1)\right] = {\mathbb E}
\left[\mathbb E \left[ h (X_1, X_2,X_3 ) | X_1 \right] \right] =
\mathbb E\left[h (X_1, X_2,X_3 )\right] = F_3 (\theta),
\end{equation}
we have
\begin{equation}
\zeta_1(F) = {\mathbb V}\text{ar}  ( h_1 (X) ) > 0
\end{equation}
unless $ F $ is degenerate. The asymptotic normality of $T_n /
{n \choose 3}$, i.e. Theorem \ref{triangle_as}(b) follows
from Theorem A of \cite{Serfling} (1980), page 192. 

As an aside, we note that the method of U-statistics employed above is more 
versatile. It may be applied to configurations involving any subgraph 
composed of finitely many vertices.  

Let $\mathbb G$ be a graph on $k$ vertices ($k \leq n$) and suppose 
$\mathbb V:= \{v_1, \dotsc , v_k\}$ and $\mathbb E$ be the vertex and edge 
sets of the graph, respectively.
For a permutation $\sigma:\mathbb V \to \mathbb V$, let us 
consider the graph $\mathbb G_\sigma$ with vertex set $\mathbb V$ and edge 
set $\mathbb E_\sigma := \{ \left<\sigma(v_i), \sigma(v_j)\right>: \left<v_i, v_j\right> \in 
\mathbb E\}$. Let $l:= \# \{ \sigma: \mathbb E_\sigma = \mathbb E\}$. Thus 
$l$ counts the number of symmetries of $\mathbb G$. To illustrate this, 
 consider the graph with vertex set as $\mathbb V = \{1,2,3,4\}$ and edge 
set $\mathbb E = \{\left<1,2\right>, \left<2,3\right>, \left<3,4\right>,  \left<4,1\right>\}$. For this graph $l= 8$; 
indeed any of the eight permutations $ (1,2,3,4), (2,3,4,1), (3,4,1,2), 
(4,1,2,3), (1,4,3,2), (4,3,2,1), (3,2,1,4) $ and $  (2,1,4,3)$ result in a graph 
which is equivalent to the original graph.

Now given a graph $G$ on the vertex set $V:=\{1,2,\dotsc , n\}$ with edge 
set $E$ consider the subgraph $G(i_1, \dotsc, i_k)$ on the vertex set $V(i_1, 
\dotsc, i_k) = \{i_1, \dotsc, i_k\}$ and edge set $E(i_1, \dotsc, i_k) = \{\left<i_t, i_s\right> 
: \left<i_t, i_s\right> \in E \text{ and } 1 \leq t\neq s \leq k\}$.
Let $ {\mathbb L}_k $ be the set of all {\em {ordered}}\/  $ k $-tuples 
from $ \{1,2, \dotsc, n \} $ and
$ f : {\mathbb L}_k \to \{ 0,1 \}$ be defined 
as follows: 
\begin{equation}
f(i_1,i_2, \dotsc, i_k) = \begin{cases} 1 & \text{ if } \left<i_s, i_t\right> \in  
E(i_1, i_2, \dotsc, i_k) \text{ for all } \left<v_s, v_t\right> \in \mathbb E  \\
0 & \text{ otherwise.} \end{cases}
\end{equation}
Thus, the number
of subgraphs in $G$ isomorphic to $\mathbb G$ is given by $l^{-1}
T_n ({\mathbb G})$, where 
\begin{equation}
\label{def:num_sub}
T_n ({\mathbb G}) := \sum_{ (i_1,i_2, \dotsc, i_k) \in {\mathbb L}_k } 
f(i_1,i_2, \dotsc, i_k) . 
\end{equation}

For the graph $G_\theta$, the vertices  $ i_s $ and $ i_t $ are connected 
by an edge  if and only if $ X_{i_s} + X_{i_t} > \theta $.  Thus the function 
$f$ in the definition (\ref{def:num_sub}) may be replaced by the  kernel 
function $ h : {\mathbb R}^k \to \{ 0,1\}$ defined by
\begin{equation}
h(X_{i_1}, X_{i_2}, \dotsc, X_{i_k})  = \begin{cases} 1 
& \text{ if } X_{i_s} + X_{i_t} > \theta \text{ for all } \left<v_s, v_t\right> \in \mathbb E  \\
0 & \text{ otherwise,} \end{cases}
\end{equation}
and we have
\begin{equation}
T_n ({\mathbb G}) = \sum_{ (i_1, i_2, \dotsc, i_k) \in {\mathbb L}_k } h 
(X_{i_1}, X_{i_2}, \dotsc, X_{i_k} ).
\end{equation}
This function $ h $ need not be symmetric in its argument. Even if it 
were so, since the sum is over all ordered $k$-tuples $(i_1, \dotsc, i_k)$ and 
not on  $k$-tuples $(i_1, \dotsc, i_k)$ such that $i_1 < i_2<
\dotsi <i_k$, 
$T_n$ need not be a U-statistic.

To overcome these,
 we first  define a symmetrized version of the kernel $h$. Let,
\begin{equation}
\label{def:sym_h}
 h_{\text{sym}} (x_1,x_2, \dotsc,x_k) := \frac{1}{k!}
 \sum_{\tilde{\sigma}
\in P(k) } h 
(x_{\tilde{\sigma}(1)}, x_{\tilde{\sigma}(2)}, \dotsc, x_{\tilde{\sigma}(k)} ) 
\end{equation}
where $\tilde{\sigma}: \{ x_1,x_2,\dotsc,x_k\} \to \{ x_1,x_2,\dotsc,x_k\} $ is a permutation,
and $P(k)$ is collection of all such permutations. We note that, even if 
$x_1 = x_2$, we take $(x_2, x_1, x_3, \dotsc , x_k)$ to be a permutation of 
$(x_1, x_2, x_3, \dotsc , x_k)$. Clearly, $ h_{\text{sym}} $ is a 
symmetric function of $ x_1, x_2, \dotsc, x_k $, and we obtain
\begin{equation}
T_n (\mathbb G) = \sum_{ (i_1, i_2, \dotsc, i_k) \in {\mathbb L}_k } h 
(x_{i_1}, x_{i_2}, \dotsc, x_{i_k} ) = \sum_{ (i_1, i_2, \dotsc, i_k) 
\in {\mathbb L}_k } h_{\text{sym}} (x_{i_1}, x_{i_2}, \dotsc, x_{i_k}
).
\end{equation}
The statistic $T_n(\mathbb G)$ is {\it not}\/ a U-statistic. Thus we
consider the statistic
\begin{eqnarray}
T'_n(\mathbb G) &=& \sum_{i_1=1}^n \sum_{i_2=1}^n \dotsi \sum_{i_k=1}^n 
  h_{\text{sym}}   (x_{i_1}, x_{i_2}, \dotsc, x_{i_k} )\nonumber\\
&=& T_n + \sum_{\begin{subarray}{1} 
{1 \leq i_1, i_2, \dotsc, i_k \leq n}\\{\text{not all distinct}}
\end{subarray}}
  h_{\text{sym}} (x_{i_1}, x_{i_2}, \dotsc, x_{i_k} )\nonumber\\
& = & T_n(\mathbb G) + R_n(\mathbb G) \text{ (say).}
\end{eqnarray}
Since $T'_n(\mathbb G)$ is a von Mises' statistic (see \cite{Koroljuk}, 
page 39), the asymptotic results about $ T'_n (
{\mathbb G}) $ can be read off from the results about the von Mises' statistics
with the kernel function $ h_{\text{sym}} $. 
To relate the statistics $T_n(\mathbb G)$ and $T'_n(\mathbb G)$ we 
observe  that the number of terms in the sum defining $R_n$ is of the 
order of $n^{k-1}$;  thus noting that $h_{\text{sym}} \leq 1$, we have 
$R_n =O(n^{k-1})$ as $n \to \infty$.
Let 
\begin{equation}
F({\mathbb G}) := {\mathbb E}\left[ h(X_{1}, X_{2}, \dotsc,
  X_{k})\right].
\end{equation}
Then, by the i.i.d. nature of $ \{ X_i : i \geq 1 \} $, we have $ {\mathbb E} 
\left[h_{\text{sym}} (X_{1}, X_{2}, \dotsc, X_{k})\right]= F({\mathbb G}) $.
\begin{theo}
\label{VE_as}
As $ n \to \infty $,
\begin{equation}
\frac{ T_n ({\mathbb G})}{{ n^k}} \to   F({\mathbb G}) 
\text{ almost surely. }
\end{equation}
\end{theo}

\noindent{\bf Proof:} 
From Theorem 3.3.1 of \cite{Koroljuk}, 
page 102, we have  $\frac{T'_n ({\mathbb G})}{{ n^k}} \to F({\mathbb G}) 
$ almost surely.  Our observation that $R_n =O(n^{k-1})$ as $n \to 
\infty$ completes the proof of the theorem.

To obtain the central limit theorem, as in (\ref{assocfn}) let 
\begin{equation}
h_1 (x) 
:= {\mathbb E}\left[ h_{\text{sym}} (x, X_2, X_3, \dotsc, X_k)\right]
\end{equation}
and $ \zeta_1({\mathbb G}) = \text{ Var} (h_1 (X) ) $ where $X$ is an 
independent random variable identical in distribution to $X_1$. Then, from 
Theorem 4.2.5 (\cite{Koroljuk}, page 134) we have the central 
limit theorem for $T'_n(\mathbb G)$. Now $R_n(\mathbb G)/n^k \to 0$ in 
probability as $n \to \infty$. Thus we obtain  
\begin{theo}
\label{VE_clt}
As $ n \to \infty $, 
\begin{equation}
\sqrt{n} \Bigl[ \frac{T_n ({\mathbb G})}{{ n^k}} - F({\mathbb G}) \Bigr]
 \Rightarrow \sqrt{k \zeta_1 ({\mathbb G})} Z.
\end{equation}
\end{theo}

\section{Local Properties} 
In this section we study $T_n(1)$ as defined in equation
(\ref{def:loc_tringles}).

For fixed $x \in {\mathbb R}$,
the kernel $ h (x,x_1,x_2) $ as defined in equation (\ref{hdef-rr})
is a symmetric function of $ x_1 $ and $ x_2 $. Define
a U-statistic based on the kernel $ h (x,\cdot,\cdot)$ by
\begin{equation}
T_n (1;x) : = \sum_{ 2 \leq i \neq j \leq n+1 } h(x,X_i,X_j).
\end{equation}
We have by the strong law for U-statistics (Theorem A, \cite{Serfling}
page 190)
\begin{equation}
\frac{T_n (1;x)}{ {n\choose2}} \to {\mathbb E}\left[h(x,X_i,X_j) \right]
\quad \text{almost surely, as}\quad n \to \infty.
\end{equation}
The random variable $ T_n (1;x) $ may be easily identified as the number of triangles
of $G_\theta$ with a fixed vertex $ 1 $ and $ X_1 = x $. Formally,
we may write, 
for any $ t \in 
{\mathbb R} $
\begin{eqnarray}
\label{thm6res}
{\mathbb E} \left[ \exp\left( i t \frac{T_n (1)}{{ n\choose 2}} \right)\right]
& = &\int_{\mathbb R} F(dx)   {\mathbb E}\left[ \exp \left( i t \frac{T_n(1)}
{{ n\choose 2}} \right)\bigg| X_1 = x\right] \nonumber\\
& \to & \int_{\mathbb R} F(dx) \exp 
\left( i t {\mathbb E} \left[ h (x,X_2,X_3) \right]\right)
\end{eqnarray}
as $n \to \infty$, because $ T_n (1;x)/{ n\choose 2} \to {\mathbb E} 
\left[h (x,X_2,X_3)\right]$ almost surely and hence also in distribution. The limit is 
justified by the usual application of dominated convergence theorem
since the integrand, being 
a characteristic function, is bounded by $ 1 $. 
The right hand side of (\ref{thm6res}) is the characteristic function of the 
random variable $ \int_{\mathbb R} \int_{\mathbb R} F(dx_1) F(dx_2) 
h(X,x_1,x_2)$. This proves the
Theorem \ref{thm:local_triangles}.   

\vskip .5em

As in the previous section we may generalize Theorem \ref{thm:local_triangles} for the subgraph $\mathbb G$ defined earlier.
 
\section{Spatial model}
Consider the Poisson spatial model as elaborated in Section 1. We
first thin the underlying Poisson process. Fix $x\in
{\mathbb R}$. For $i \geq 1$, each point
$ \xi_i$ of the original process is included in the thinned process
with a probability $ f (|\xi_i|;x)$ independently of the other points.
The thinned process is an inhomogeneous Poisson process with intensity
function $g(y)$ for $y\in{\mathbb R}$ given by 
$g(y)= \lambda f(|y|;x)$, where $f(r;x)$ is as defined in
equation~(\ref{r-fdef}). For $\Delta_r$ as defined
in equation~(\ref{def:deg0n}), we have
\begin{prop}
\label{distrn}
The conditional distribution of $\Delta_r$ given that
$ X_0 = x $, is Poisson with parameter $ \lambda_r
(x) $ where 
\begin{equation}
\lambda_r (x) : = \lambda
\int_{ |\tilde{r}| \leq r } f(\tilde{r};x) d\tilde{r}.
\end{equation}
\end{prop}

We first prove Theorem \ref{finitecase} where $ C_r (x) \to C(x) <
 \infty$ (see Section 1 for the relevant definitions). Because $ t \to
 \int_{\mathbb R} F(dx) \exp \bigl( - \lambda c_d C(x) ( 1 - \exp (it)
 ) \bigr) $ is indeed a characteristic function, it is enough to prove
 that the characteristic function of $\Delta_r$ converges to the above
 quantity.
By Proposition \ref{distrn}, the conditional distribution of $\Delta_r$ given
$ X_0 = x $ is Poisson with parameter $ \int_{|\tilde{r}|\leq r}
\lambda f (\tilde{r},x) d\tilde{r}
=  \lambda c_d  \int_0^r  \tilde{r}^{d-1} [ 1 - F( \theta 
\tilde{r}^{\beta} - x)] d\tilde{r}
= \lambda c_d C_r(x)$ where $ c_d $ is the volume of the $(d-1)$-dimensional
unit sphere. 
Therefore, we have
\begin{eqnarray}
 \phi_{\Delta_r} (t) & := & {\mathbb E} \left[ \exp (it
 \Delta_r)\right]\nonumber\\
& = & \int_{\mathbb R} F(dx)  {\mathbb E} \left[
\exp (it \Delta_r) | X_0 = x\right]\nonumber\\
& = &  \int_{\mathbb R} F(dx) \exp \bigl( - \lambda c_d C_r (x) ( 1 - 
\exp (it) ) \bigr).
\end{eqnarray}
Now, since $C_r (x) \to C(x)$, 
the usual dominated convergence theorem assures
\begin{equation}
\phi_{\Delta_r} (t) \to \phi_{\Delta} (t)  \text{ as } r \to \infty.
\end{equation}
This completes the proof of Theorem \ref{finitecase}. 

\vskip 1em
\noindent{\bf Remark}
Assuming $ \theta > 0 $,  we may re-write 
$ C(x) $ in the following way:
\begin{eqnarray}
\lefteqn{ \int_0^{\infty} r^{d-1} [ 1 - F( \theta r^{\beta} -x ) ] dr
}\nonumber\\
& = & \int_0^{\infty} r^{d-1} \int_{-\infty}^{ \infty} 1 \bigl( \tilde{x} > \theta 
r^\beta -x \bigr) F (d\tilde{x}) dr\nonumber\\
& = & \int_{-\infty}^{ \infty} \int_0^{\infty} r^{d-1} 1 \bigl( r < 
\max \{ 0, [(\tilde{x}+x)/\theta]^{1/\beta}\} \bigr) dr F(d\tilde{x})
\nonumber\\
& = & \frac{1}{d}  \int_{-\infty}^{ \infty}  \bigl[ \max \{ 0, [(\tilde{x}+x)/\theta
]^{d/\beta} \} \bigr]  F(d\tilde{x})\nonumber\\
& = & \frac{1}{d \theta^{d/\beta} }  \int_{-\infty}^{ \infty}  \bigl[
  \max \{ 0, (\tilde{x}+x)^{d/\beta} \} \bigr]
F(d\tilde{x})\nonumber\\
& = &   \frac{1}{d \theta^{d/\beta} }  {\mathbb E} \bigl[ \max\{ 0, (X_0 
+ x)^{d/\beta} \} \bigr].
\end{eqnarray}
Thus $ C(x) < \infty $ if  $ {\mathbb E} [ |X_0| ^{d/\beta} ] < \infty $.

\vskip 1em
To show Theorem \ref{clt}, it is enough to prove that the characteristic function
of $ (\Delta_r - \lambda c_d C_r)/\sqrt{ \lambda c_d C_r}$ converges to 
the product of the characteristic functions of a standard normal
random variable and $\sqrt{\lambda c_d} g(X_0)$. 
\begin{eqnarray}
& & {\mathbb E} \left[
\exp\left( it \frac{\Delta_r - \lambda c_d C_r}{\sqrt{\lambda c_d 
C_r}} \right)\right]\nonumber\\
& = &  \int_{\mathbb R} F(dx) {\mathbb E} \left[ \exp \left(
it \frac{\Delta_r - 
\lambda c_d C_r}{\sqrt{\lambda c_d C_r}} \right) \bigg| X_0 = x \right]\nonumber\\
& = &  \int_{\mathbb R} F(dx) \exp ( -it \sqrt{ \lambda c_d C_r} ) 
{\mathbb E} \left[ \exp \left( i \frac{ t }{\sqrt{\lambda c_d C_r} }
\Delta_r \right)  
\bigg| X_0 = x \right]\nonumber\\
& = & \int_{\mathbb R} F(dx) \exp \left( -it \sqrt{ \lambda c_d C_r}  
-  \lambda c_d C_r(x) ( 1 - \exp( it / \sqrt{ \lambda c_d C_r} )) \right)
\end{eqnarray} 
using the fact that the conditional distribution of $\Delta_r$ given
 $ X_0 = x $
follows a Poisson distribution with parameter $  \lambda c_d C_r(x)$. 
Consider the  
logarithm of the integrand:
\begin{eqnarray}
\lefteqn{ -it \sqrt{ \lambda c_d C_r} - \lambda c_d C_r (x) \left[ 1 -
    \exp \left( \frac{it}{ \sqrt{ \lambda c_d C_r} } \right)
    \right]}\nonumber\\
& = & -it \sqrt{ \lambda c_d C_r} - \lambda c_d C_r(x) \left[ - \frac{
    it }{ \sqrt{ \lambda c_d C_r}} - \frac{1}{2} \left( \frac{it}{
    \sqrt{ \lambda c_d C_r}} \right)^2 + o \left( \left(\frac{it}{
    \sqrt{ \lambda c_d C_r} }\right)^2 \right) \right]\nonumber\\
&=& it \sqrt{ \lambda c_d}\Bigl[ \frac{ C_r(x) - C_r}{ \sqrt{C_r}}\Bigr] 
- \frac{t^2}{2} \frac{ C_n (x) }{ C_r} - \lambda c_d C_r(x)  \times 
o\Bigl( \frac{1}{ C_r}\Bigr).
\end{eqnarray}
Under our condition, the first term converges to $ it \sqrt{ \lambda
c_d}g(x) $. The condition also implies that $C_r(x)/C_r\to
1$ as $r\to \infty$, and thus the second term converges to $-t^2/2$.
The third term can be written as $ \bigl[ C_r(x) / C_r\bigr] \times
\bigl[ C_r o(1/C_r) \bigr] \to 0 $ as $r\to\infty$.  Applying the
dominated convergence theorem, we now obtain that
\begin{equation}
{\mathbb E} \left[ \exp\left(
it \frac{ \Delta_r - \lambda c_d C_r}{\sqrt{\lambda c_d 
C_r}} \right)\right]
\to
\exp\left(-\frac{t^2}{2}\right)
\int_{\mathbb R} F(dx)\exp(it \sqrt{\lambda c_d} g(x)).
\end{equation}
Now note that $\exp(-t^2/2)$ is the characteristic function of a standard normal 
random variable  and 
$\int_{\mathbb R} F(dx)\exp(it \sqrt{\lambda c_d} g(x))$ is the characteristic 
function of $\sqrt{\lambda c_d} g(X_0)$. This completes the proof of 
Theorem \ref{clt}. 

\section*{Acknowledgments}

We thank Hiroyoshi Miwa for helpful discussion.
Naoki Masuda is supported by Special
Postdoctoral Researchers Program of RIKEN, and
Rahul Roy is supported by DST Grant MS/152/01.

\end{document}